\documentclass[English,11pt]{amsart}

\usepackage{amssymb}
\usepackage{amsmath}
\usepackage{amsfonts}
\usepackage{tabularx}
\usepackage{frcursive}
\usepackage{enumerate}
\usepackage[usenames, dvipsnames]{xcolor}
\usepackage{pgf,tikz}
\usetikzlibrary{arrows}
\usepackage{xcolor}
\usepackage{color}
\usepackage{hyperref}
\hypersetup{
    colorlinks=true,
    linkcolor=red,
    filecolor=magenta,      
    urlcolor=cyan,
}
\newtheorem{thm}{Theorem}
\newtheorem{prop}[thm]{Proposition}
\newtheorem{lem}[thm]{Lemma}

\newcommand{\IET}{\mathrm{IET}} 
\newcommand{\AIET}{\mathrm{AIET}} 
\newtheorem{definition}[thm]{Definition}
\newtheorem{example}[thm]{Example}
\newtheorem{remark}[thm]{Remark}

\newcommand{\R}{\mathbb{R}}
\newcommand{\N}{\mathbb{N}}
\newcommand{\Z}{\mathbb{Z}}
\newcommand{\Q}{\mathbb{Q}}

\usepackage[T1]{fontenc}

\begin{document}

\title{Relations with a fixed interval exchange transformation
}


\author[1]{Magali Jay}\thanks{M. \textsc{Jay} \\ Aix Marseille Université, Institut de Mathématiques de Marseille, I2M - UMR 7373, 13453 \\
Marseille, France \\
E-mail: magali.jay@univ-amu.fr}


\maketitle

\begin{abstract}
We study the group of all interval exchange transformations (IETs). Katok asked whether it contains a free subgroup. We show that for every IET $S$, there exists a dense open set $\Omega(S)$ of admissible IETs 
such that the group generated by $S$ and any $T\in\Omega(S)$ is not free of rank 2. This extends a result by Dahmani, Fujiwara and Guirardel in \cite{DFG}: the group generated by a generic pair of elements of IET([0;1)) is not free (assuming a suitable condition on the underlying permutation).

\end{abstract}

\section{Introduction}
\label{intro}

In the Oberwolfach Conference \textit{Geometric Group Theory, Hyperbolic Dynamics and
Symplectic Geometry}, which took place on September 2008, Katok \cite{KatokQuestion} posed the question whether the group of interval exchange transformations (IETs) contains a rank 2 free subgroup.

The history of IETs begun with the first return map of flows, introduced by Poincaré, although they were not formalised at this time. In the 60's, russian mathematicians worked on them in the term of differential forms. Sataev \cite{Sataev} showed that there exists a minimal IET with $g$ invariant ergodic mesures, for all integer $g$.

In 1967, Katok and Stepin \cite{KatokStepin} formally introduced 3-IETs. Keane \cite{Keane75} studied the IETs with an arbitrary number of intervals, which have subsenquently been widely investigated since 1980.
While IETs with reducible underlying permutation are of course not minimal, Keynes-Newton \cite{KeynesNewton}, Keane \cite{Keane77} and later Coffrey \cite{Coffrey}, gave examples of minimal but non uniquely ergodic IETs.
In the early days of the study of IETs, Veech \cite{Veech} had also given examples of minimal and non uniquely ergodic dynamical systems, which where related to IETs, although he did not write his article in the terms of IETs.

IETs arise in many areas of mathematics such as dynamical systems, polygonal billiards, geometry and flows on flat surfaces. See \cite{Viana} for a survey on the subject. They are well understood individually as maps, but the structure of their group still gives rise to open questions, such as that of Katok. Currently, we know certain results e.g. about subgroups generated by torsion elements and rotations \cite{Boshernitzan}. Its commutator group is a proper normal subgroup (independently shown by Arnoux \cite{Arnoux_eif} and Sah \cite{Sah}). Therefore the group of interval exchange transformations is not simple, in contrast to the group of interval exchange transformations with flips, (see part III.1 of \cite{ArnouxFIET}). Moreover, the commutator is exactly the subgroup generated by torsion elements \cite{Vorobets}.
 
The main result of this paper gives a partial answer to Katok's question:

\begin{thm} \label{thJ} Let $S$ be any IET on $[0;1)$.
There exists a dense open set $\Omega_{a}(S) \subset \IET_{a}([0;1))$ such that for every $T \in \Omega_{a}(S)$, $\langle S,T \rangle$ is not free of rank 2.
\end{thm}

The set $\IET_a([0;1))$ is the subset of admissible IETs (see Definition \ref{def_admissible}, this is the same definition as in \cite{DFG}). This theorem extends the following result, proved in \cite{DFG}.
\begin{thm}[Dahmani-Fujiwara-Guirardel]\label{thDFG}
There exists a dense open set $\Omega \subset \IET([0;1)) \times \IET_{a}([0;1))$ such that for every $(S,T) \in \Omega$, $\langle S,T \rangle$ is not free of rank 2.
\end{thm}

We also present some computations to determine whether self-similar IETs, such as the Arnoux-Yoccoz example (see \cite{AY} and \cite{Arnoux}), can be in some $\Omega_{a}(S)$.

Note that Katok's question is easily answered in the group of affine interval exchange transformations (AIETs), see Definition \ref{def_AIET}.

\begin{prop} \label{free_AIET}
There exist $f,g \in \AIET([0;1))$ such that $\langle f,g \rangle$ is a free group of rank 2.
\end{prop}

\noindent
\textbf{Acknowledgements:} 
This article comes from my Master's thesis supervised by Erwan Lanneau at the Institute Fourier of Grenoble Alpes University. I thank Erwan Lanneau for the topic, his advice and discussions.
I am also thankful to Pierre Arnoux and Vincent Delecroix for discussions about this subject in Autrans, and to Adrien Boulanger, who pointed out the free group in the group of AIETs.
I thank Diana Davis for her help with English language.
Last but not least, I thank the anonymous reviewer for their careful reading and helpful comments.





\section{Outline of the proof}

Before going into detail we summarise the proof of Theorem \ref{thJ}:
\begin{enumerate}
\item
(\emph{Basic relations}).
Let $T_0$ be a $q$-rational IET (meaning that all its discontinuities are at rational points with denominator $q$), and two obvious relations between $S$ and $T_0$. Say $r_1(S,T_0) = \mathrm{id} = r_2(S,T_0)$.

\item
(\emph{Small pertubations}).
Choose any $T$ close to $T_0$. Then $r_1(S,T)$ and $r_2(S,T)$ induce translations with small translation length on every interval of continuity that is not too close to some finite set of points $X_q(S)$ (see Definition \ref{defXYZ}).

\item
(\emph{IET with small support}).
The commutator $U = [r_1(S,T),r_2(S,T)]$ has a small support, located in a neighbourhood of $X_q(S)$.

\item
(\emph{Drifting the support}).
With an additional condition on $T_0$, some power $k$ of $T$ is such that $T^k (supp(U)) \cap supp(U) = \emptyset$. We say that $T$ \emph{drifts} the support of $U$.

\item
(\emph{Relation}).
Then $U$ commutes with $T^{k}UT^{-k}$, which gives the nontrivial relation $[U,T^{k}UT^{-k}] = \mathrm{id}$.

\end{enumerate}

Parts (1) through (3) are done in section \ref{section_small_supp}, point (4) in section \ref{section_drift_supp} and part (5) in section \ref{ccl}.

The main differences between this work and the proof of Theorem \ref{thDFG} by Dahmani, Fujiwara and Guirardel lies in (4). If $S$ is rational, we can take $X_q = \frac{1}{q}\N$ (as implicitly in \cite{DFG}) but whenever $S$ is not rational, this cannot be done anymore. This is why our definition of $X_q(S)$ generalises that of \cite{DFG} and coincides with it if and only if $S$ is rational.

Before stating all lemmas that enable us to prove Theorem \ref{thJ}, we need to recall some classical notation and introduce some others.

\section{Background for the proof}
\subsection{Classical notation}
We do not recall the definition of \emph{IET}, nor related definitions of \emph{interval of continuity}, \emph{length of translation}, \emph{underlying permutation}. See \cite{Viana} for a survey; we use the same notation. First of all we define a few linear maps, useful when speaking about 
IETs.

\begin{definition} Let $T$ be an IET on $[0;1)$. We denote by $\lambda(T)=(l_1,...,l_n)$ the set of lengths of its intervals of continuity, by $\beta(T)=(b_1,...,b_{n-1})$ its points of discontinuity, by $\pi(T)$ its underlying permutation and by $\omega(T)=(t_1,...,t_n)$ its translation lengths.
\end{definition}

If needed, we write $b_0 =0$ and $b_n =1$.
To summarise, here are the four maps we will use in this article:
$$
\begin{array}{rccl}
\lambda: & \IET([0;1)) & \twoheadrightarrow & \left\{ (a_1,...,a_n) \ \left| 
\begin{array}{ll}
\ n\in\N, & \forall i \in \{1,...,n\}, \, a_i \in \R_+^*, \\
 & \sum_{i=1}^{n} a_i = 1
\end{array}
\right.
\right\} \\
 & T & \mapsto & \lambda(T) = (l_1,...,l_n) \\ [0.5cm]
\beta: & \IET([0;1)) & \twoheadrightarrow & \{(a_1,...,a_{n-1}) \ | \ n\in\N, \, 0 < a_1 < ... < a_{n-1} < 1 \} \\
 & T & \mapsto & \beta(T) = (b_1,...,b_{n-1}) \\ [0.5cm]
\pi: & \IET([0;1)) & \twoheadrightarrow & \left\{ \sigma \in \mathfrak{S}_n \ \left| \ \begin{array}{l}
n\in\N, \\ 
\forall i \in \{1,...,n-1\}, \sigma(i+1) \neq \sigma(i) + 1
\end{array} \right. \right\} \\
& T & \mapsto & \pi(T) \\ [0.5cm]
\omega: & \IET([0;1)) & \rightarrow & \{(a_1,...,a_n) \ | \ n\in\N, \,\forall i, a_i \in (-1;1)  \} \\
 & T & \mapsto & \omega(T) = (t_1,...,t_n)
\end{array}
$$


We will denote by $\Delta(T)$ the set of discontinuity points of an IET $T$. If $X$ is a set, let $\mathcal{N}_\epsilon(X) = \{ y \, | \, d(y,X) < \epsilon \}$ be the open $\epsilon$-neighbourhood of $X$. 

We denote by $\IET_\sigma$ the set $\pi^{-1}(\sigma)$ of IETs having $\sigma$ as their underlying permutation.
Finaly, we equip $\IET([0;1))$ with the following distance:
$$
\begin{array}{l c c l}
d :& \IET([0;1))\times \IET([0;1)) & \longrightarrow &\R_+\cup\{\infty\} \\
 & (S,T) & \longmapsto &
\left\{
\begin{array}{l l}
 || \lambda(S) - \lambda(T)||_1 & \text{if } \pi(S)= \pi(T)\\
 \infty & \text{otherwise}
\end{array}
\right.
\end{array}
.
$$

As in \cite{DFG}, we define the term \emph{admissible}:

\begin{definition} \label{def_admissible}
A permutation  $\sigma	\in \mathfrak{S}_n$ is said to be \emph{non admissible} when there exists  $k \in \{1,..., n\}$, such that $\sigma(k) =k$ and $\sigma(\{1,...,k\}) = \{1,...,k\}$.  
An IET $T$ is \emph{admissible} when $\sigma	= \pi(T) \in \mathfrak{S}_n$ is an admissible permutation.

\noindent
We denote by $\IET_{a}$ the set of admissible IETs on $[0;1)$.
\end{definition}

\subsection{Specific notation}

\begin{definition} \label{defXYZ} Let $S$ be any IET. Let $q$ be a positive integer.
Let us define
$$X_q(S)  = \Delta(S^{-1}) \cup S\left(\frac{1}{q} \N \right) \cup \frac{1}{q}\N$$
and
$$ Y_q(S) = \pi_q (X_q(S)) $$
 where $\pi_q : [0;1] \rightarrow \left[0;\frac{1}{q}\right) $ is the canonical projection modulo $\frac{1}{q}$.

Also define
$$Z_q(S) = Y_q(S)\sqcup \left( \frac{1}{q} + Y_q(S) \right) \sqcup \left( \frac{2}{q} + Y_q(S) \right) \sqcup \dots \sqcup \left( \frac{q-1}{q} + Y_q(S) \right) \sqcup \{1\} .$$

Denote by 
$$\alpha_q(S) = \underset{\text{C connected component of }[0 ; \frac{1}{q}) \, \setminus \, Y_q(S)}{\text{min}} \text{diam}(C)$$
the length of the smallest interval of $\left[0 ; \frac{1}{q}\right) \, \setminus \, Y_q (S)$.

Finally define
$$\mathcal{U}_{\epsilon}^q = \{ R \in \IET([0;1)) \, | \, \alpha_q(R) > \epsilon \}.$$
\end{definition}

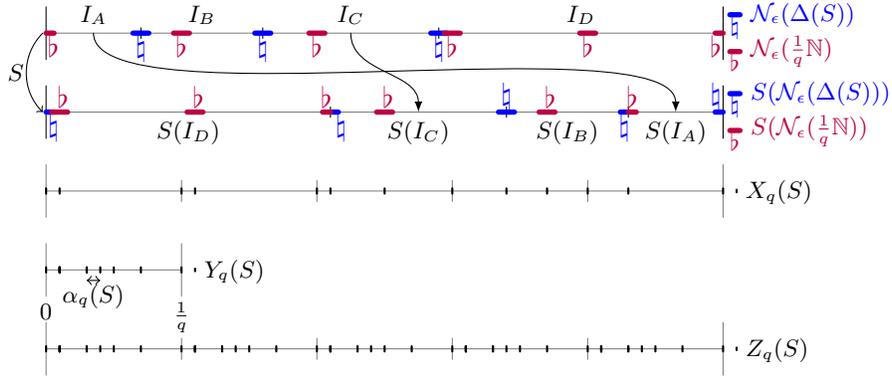
\begin{figure}[h!]
\center
\begin{tikzpicture}[line cap=round,line join=round,x=0.90cm,y=0.7cm]
\def \s{1.5} 
\def \h{1.5} 
\def \t{4}   
\def \A{1.4} 
\def \B{1.8} 
\def \C{2.6} 
\def \D{4.2} 
\def \L{\A+\B+\C+\D} 
\def \q{2} 
\def \eps{3} 
\def \colS{blue}  
\def \colN{purple}  

\draw[color=gray] (\s,0) -- (\s+\L,0);
\foreach \x in {\A,\A+\B,\A+\B+\C}
{
\draw[shift={(\x+\s,0)},color=black] (0pt,2pt) -- (0pt,-2pt);
\draw[shift={(\x+\s,0)},color=\colS,line width=2pt] (-\eps pt,0pt) -- (\eps pt,0pt);
\draw[shift={(\x+\s,-0.3)},color=\colS] node {$\natural$};
}
\draw (\s,-0.5) -- (\s,0.5);
\draw (\s+\L,-0.5) -- (\s+\L,0.5);

\draw[color=gray] (\s,-\h) -- (\s+\L,-\h);
\foreach \x in {\D,\D+\C,\D+\C+\B}
\draw[shift={(\x+\s,-\h)},color=black] (0pt,2pt) -- (0pt,-2pt);

\draw (\s,-\h-0.5) -- (\s,-\h+0.5);
\draw (\s+\L,-\h-0.5) -- (\s+\L,-\h+0.5);

\foreach \d in {2*\h,4*\h}
{
\draw[shift={(\s,-\d)},color=gray] (0,0) -- (\L,0);
\draw[shift={(\s,-\d)},color=gray] (0,-0.5) -- (0,0.5);
\draw[shift={(\s,-\d)},color=gray] (\L,-0.5) -- (\L,0.5);
\foreach \x in {2,4,6,8}
\draw[shift={(\s+\x,-\d)},color=gray] (0pt,4pt) -- (0pt,-4pt);
}

\draw[shift={(\s,-3*\h)},color=gray] (0,0) -- (\q,0);
\draw[shift={(\s,-3*\h)},color=gray] (0,-0.5) -- (0,0.5);
\draw[shift={(\s,-3*\h)},color=gray] (\q,-0.5) -- (\q,0.5);

\foreach \x in {2,4,6,8}
{
\draw[shift={(\x+\s,0)},color=\colN,line width=2pt] (-\eps pt,0pt) -- (\eps pt,0pt);
\draw[shift={(\x+\s,-0.25)},color=\colN] node {$\flat$};
}
\draw[shift={(\s,0)},color=\colN,line width=2pt] (0pt,0pt) -- (\eps pt,0pt);
\draw[shift={(\s+0.05,-0.25)},color=\colN] (\eps/2 pt,0) node {$\flat$};
\draw[shift={(\L+\s,0)},color=\colN,line width=2pt] (-\eps pt,0pt) -- (0pt,0pt);
\draw[shift={(\L+\s-0.05,-0.25)},color=\colN] (-\eps/2 pt,0) node {$\flat$};

\foreach \d in {0}
{
\foreach \x in {\L,\D+\C+\B,\D+\C}
\draw[shift={(\x+\s,\d-\h)},color=\colS,line width=2pt] (-\eps pt,0pt) -- (0pt,0pt);

\foreach \x in {\D+\C,\D,0}
\draw[shift={(\x+\s,\d-\h)},color=\colS,line width=2pt] (0pt,0pt) -- (\eps pt,0pt);

\draw[shift={(\D+\s+0.1,\d-\h-0.3)},color=\colS] (\eps/2 pt,0) node {$\natural$};
\draw[shift={(\D+\C+\s,\d-\h+0.3)},color=\colS] (0,0) node {$\natural$};
\draw[shift={(\s+0.05,\d-\h-0.3)},color=\colS] (\eps/2 pt,0) node {$\natural$};
\draw[shift={(\D+\C+\B+\s,\d-\h-0.3)},color=\colS] (-\eps/2 pt,0) node {$\natural$};
\draw[shift={(\D+\C+\B+\A+\s-0.05,\d-\h+0.3)},color=\colS] (-\eps/2 pt,0) node {$\natural$};

\foreach \x in {\q-\A+\D+\C, 2*\q-\A-\B+\D, 3*\q-\A-\B-\C,4*\q-\A-\B-\C}
{
\draw[shift={(\x+\s,\d-\h)},color=\colN,line width=2pt] (-\eps pt,0pt) -- (\eps pt,0pt);
\draw[shift={(\x+\s+0.05,\d-\h+0.3)},color=\colN] (0,0pt) node {$\flat$};
}

\draw[shift={(\s+\D+\C+\B,\d-\h)},color=\colN,line width=2pt] (0pt,0pt) -- (\eps pt,0pt);
\draw[shift={(\s+\D+\C+\B,\d-\h+0.3)},color=\colN] (\eps/2 pt,0pt) node {$\flat$};
\draw[shift={(\D+\s,\d-\h)},color=\colN,line width=2pt] (-\eps pt,0pt) -- (0pt,0pt);
\draw[shift={(\D+\s,\d-\h+0.3)},color=\colN] (-\eps/2 pt,0pt) node {$\flat$};
}

\foreach \x in {0,\D,\D+\C,\D+\C+\B,\L,\q-\A+\D+\C, 2*\q-\A-\B+\D, 3*\q-\A-\B-\C,4*\q-\A-\B-\C,2,4,6,8}
\draw[shift={(\x+\s,-2*\h)},thick] (0,0.05) -- (0,-0.05);

\foreach \x in {0,
\D+\C+\B-4*\q,
\D+\C-3*\q,
\D-2*\q,
\q-\A+\D+\C-3*\q, 
2*\q-\A-\B+\D-2*\q, 
3*\q-\A-\B-\C,
4*\q-\A-\B-\C-\q,
\D+\C+\B-4*\q,
\D-2*\q,2}
\draw[shift={(+\x+\s,-3*\h)},thick] (0,0.05) -- (0,-0.05);

\draw [<->] (\s+\D+\C-3*\q, -3*\h-0.2) -- (\s+\D+\C+\B-4*\q,-3*\h-0.2);
{\footnotesize
\draw (\s+\D+\C+\B/2-3.5*\q,-3*\h-0.5) node {$\alpha_q(S)$};
\draw (\s,-3*\h-0.8) node {$0$};
\draw (\s+\q,-3*\h-0.9) node {$\frac{1}{q}$};
}

\foreach \z in {0, \q, 2*\q, 3*\q, 4*\q}
{\foreach \x in {0,
\D+\C+\B-4*\q,
\D+\C-3*\q,
\D-2*\q,
\q-\A+\D+\C-3*\q, 
2*\q-\A-\B+\D-2*\q, 
3*\q-\A-\B-\C,
4*\q-\A-\B-\C-\q,
\D+\C+\B-4*\q,
\D-2*\q}
\draw[shift={(\z+\x+\s,-4*\h)},thick] (0,0.05) -- (0,-0.05);
}
\draw[shift={(\s+5*\q,-4*\h)},thick] (0,0.05) -- (0,-0.05);

\begin{footnotesize}
\draw (\s+\A/2,0) node (A) [above]{$I_A$};
\draw (\s+\A+\B/2,0) node (B) [above]{$I_B$};
\draw (\s+\A+\B+\C/2,0) node (C) [above]{$I_C$};
\draw (\s+\A+\B+\C+\D/2,0) node (D) [above]{$I_D$};

\draw (\s+\D+\C+\B+\A/2,-\h) node (A') [below]{$S(I_A)$};
\draw (\s+\D+\C+\B/2,-\h) node (B') [below]{$S(I_B)$};
\draw (\s+\D+\C/2,-\h) node (C') [below]{$S(I_C)$};
\draw (\s+\D/2,-\h) node (D') [below]{$S(I_D)$};

\draw[->,>=latex] (A) to[out=-90, in =90,distance=1.2cm] (A');
\draw[->,>=latex] (C) to[out=-90, in =90] (C');

\draw [->] (\s-0.05,0) to [out = -135,in =135] (\s-0.05,-\h);
\draw (\s-0.45,-\h/2) node {$S$}; 
\end{footnotesize}
\begin{footnotesize}
\draw[\colS,line width=2pt] (\s+\L+0.1,0.35) -- (\s+\L+0.25,0.35) node [right]{$\mathcal{N}_{\epsilon}(\Delta(S))$};
\draw[shift={(\s+\L+0.175,0.1)},color=\colS] node {$\natural$};
\draw[\colN,line width=2pt] (\s+\L+0.1,-0.35) -- (\s+\L+0.25,-0.35) node [right]{$\mathcal{N}_{\epsilon}(\frac{1}{q}\N)$};
\draw[shift={(\s+\L+0.175,-0.55)},color=\colN] node {$\flat$};

\draw[\colS,line width=2pt] (\s+\L+0.1,-\h+0.35) -- (\s+\L+0.25,-\h+0.35) node [right]{$S(\mathcal{N}_{\epsilon}(\Delta(S)))$};
\draw[shift={(\s+\L+0.175,-\h+0.1)},color=\colS] node {$\natural$};
\draw[\colN,line width=2pt] (\s+\L+0.1,-\h-0.35) -- (\s+\L+0.25,-\h-0.35) node [right]{$S(\mathcal{N}_{\epsilon}(\frac{1}{q}\N))$};
\draw[shift={(\s+\L+0.175,-\h-0.55)},color=\colN] node {$\flat$};

\draw[thick] (\s+\L+0.2,-2*\h+0.025) -- (\s+\L+0.2,-2*\h-0.025) node [right]{$X_q(S)$};
\draw[thick] (\s+\q+0.2,-3*\h+0.025) -- (\s+\q+0.2,-3*\h-0.025) node [right]{$Y_q(S)$};
\draw[thick] (\s+\L+0.2,-4*\h+0.025) -- (\s+\L+0.2,-4*\h-0.025) node [right]{$Z_q(S)$};
\end{footnotesize}
\end{tikzpicture}
\caption{Illustration of the notation (here, $q=5$)}
\label{fig_def_XYZ}
\end{figure}

Note that $X_q(S) \subset Z_q(S)$ and that $\alpha_q(S) \geqslant 0$ is well defined because the set $X_q(S)$ (hence $Y_q(S)$) is finite. Moreover $\alpha_q(S) = 0 $ if and only if $S$ is $q$-rational. See Figure \ref{fig_def_XYZ}. 

\begin{remark}\label{U_ngb}
Let $S$ be an IET such that for every distinct $x,y\in \Delta(S^{-1}) \cup \{0\}$, we have $x+\frac{1}{q}\N \neq y+\frac{1}{q}\N$.
If $\epsilon<\alpha_q(S)$, then the set $\mathcal{U}_{\epsilon}^q$ is a neighbourhood of $S$.
\end{remark}

Here is an example of an IET $S$ such that $\mathcal{U}_\epsilon^q (S)$ is not a neighbourhood of $S$ to show the importance of the hypothesis on $S$.

\begin{example}
Let $q\geqslant 3$ and consider $S$ defined by
$\pi(S)= \left(
\begin{array}{ccc}
A&B&C\\
B&A&C
\end{array}
\right)$ 
and $\lambda(S)=\left(\lambda_A = \frac{3}{2q}, \lambda_B = \frac{1}{q}, \lambda_C = \frac{2q-5}{2q}\right)$.  

Let $\epsilon\in]0;\frac{1}{2q})$, $\eta\in]0;\epsilon)$ and define $R$ by $\lambda(R)= \lambda(S) + (\eta,-\eta,0)$ and $\pi(R)=\pi(S)$.
See Figure \ref{counterexUalpha} (where $q=5$).

\begin{figure}[h!]
\center
\begin{tikzpicture}[line cap=round,line join=round,x=0.55cm,y=0.4cm]
\def \a{3}
\def \b{2}
\def \c{5}
\def \cola{purple}
\def \colb{orange}

\draw[\cola] (10,0) node {$[$};
\draw[\cola] (10,0) -- (10+\a,0);
\draw[\colb] (10+\a,0) node {$[$};
\draw[\colb] (10+\a,0) -- (10+\a+\b,0);
\draw (10+\a+\b,0) node {$[$};
\draw (10+\a+\b,0) -- (20,0);
\foreach \h in {0,-3}
{\foreach \x in {2,4,6,8}
\draw[shift={(10+\x,\h)}, blue, densely dotted] (0pt,5pt) -- (0pt,-5pt);}
\draw[gray] (10,-0.5) -- (10,0.5);
\draw[gray] (20,-0.5) -- (20,0.5);
\draw[\colb] (10,-3) node {$[$};
\draw[\colb] (10,-3) -- (10+\b,-3);
\draw[\cola] (10+\b,-3) node {$[$};
\draw[\cola] (10+\b,-3) -- (10+\a+\b,-3);
\draw (10+\a+\b,-3) node {$[$};
\draw (10+\a+\b,-3) -- (20,-3);
\foreach \x in {\b,\b+\a}
\draw[gray] (10,-3.5) -- (10,-2.5);
\draw[gray] (20,-3.5) -- (20,-2.5);


\begin{footnotesize}
\draw (10+\a/2,0) node (A) [above]{$I_A$};
\draw (10+\a+\b/2,0) node (B) [above]{$I_B$};
\draw (10+\a+\b+\c/2,0) node (C) [above]{$I_C$};

\draw (10+\b/2,-3) node (B') [below]{$S(I_B)$};
\draw (10+\b+\a/2,-3) node (A') [below]{$S(I_A)$};
\draw (10+\b+\a+\c/2,-3) node (C') [below]{$S(I_C)$};

\draw[blue] (12+0.15,0) node [below] {$\frac{1}{5}$};
\draw[blue] (14+0.15,0) node [below] {$\frac{2}{5}$};
\draw[blue] (16+0.15,0) node [below] {$\frac{3}{5}$};
\draw[blue] (18+0.15,0) node [below] {$\frac{4}{5}$};
\draw[blue] (12+0.15,-3) node [above] {$\frac{1}{5}$};
\draw[blue] (14+0.15,-3) node [above] {$\frac{2}{5}$};
\draw[blue] (16+0.15,-3) node [above] {$\frac{3}{5}$};
\draw[blue] (18+0.15,-3) node [above] {$\frac{4}{5}$};
\end{footnotesize}
\draw[->,>=latex] (A) to[out=-80, in =100, distance=1.05cm] (A');
\draw[->,>=latex] (B) to[out=-100, in =80] (B');
\draw[->,>=latex] (C) to[out=-90, in =90] (C');

\def \deb{23}
\def \e{3.2}
\def \f{1.8}
\def \g{5}
\def \cole{red}
\def \colf{brown}

\draw[gray] (\deb,-0.5) -- (\deb,0.5);
\draw[gray] (\deb+10,-0.5) -- (\deb+10,0.5);
\draw[\cole] (\deb,0) node {$[$};
\draw[\cole] (\deb,0) -- (\deb+\e,0);
\draw[\colf] (\deb+\e,0) node {$[$};
\draw[\colf] (\deb+\e,0) -- (\deb+\e+\f,0);
\draw (\deb+\e+\f,0) node {$[$};
\draw (\deb+\e+\f,0) -- (\deb+10,0);
\foreach \h in {0,-3}
{\foreach \x in {2,4,6,8}
\draw[shift={(\deb+\x,\h)}, blue, densely dotted] (0pt,5pt) -- (0pt,-5pt);}

\draw[gray] (\deb,-3.5) -- (\deb,-2.5);
\draw[gray] (\deb+10,-3.5) -- (\deb+10,-2.5);
\draw[\colf] (\deb,-3) node {$[$};
\draw[\colf] (\deb,-3) -- (\deb+\f,-3);
\draw[\cole] (\deb+\f,-3) node {$[$};
\draw[\cole] (\deb+\f,-3) -- (\deb+\e+\f,-3);
\draw (\deb+\e+\f,-3) node {$[$};
\draw (\deb+\e+\f,-3) -- (\deb+10,-3);


\begin{footnotesize}
\draw (\deb+\e/2,0) node (E) [above] {$I_A'$};
\draw (\deb+\e+\f/2+0.2,0) node (F) [above]{$I_B'$};
\draw (\deb+\e+\f+\g/2,0) node (G) [above]{$I_C'$};
\draw (\deb+\f/2,-3) node (F') [below]{$R(I_B')$};
\draw (\deb+\f+\e/2,-3) node (E') [below]{$R(I_A')$};
\draw (\deb+\f+\e+\g/2,-3) node (G') [below]{$R(I_C')$};

\draw[blue] (\deb+2+0.15,0) node [below] {$\frac{1}{5}$};
\draw[blue] (\deb+4-0.15,0) node [below] {$\frac{2}{5}$};
\draw[blue] (\deb+6+0.15,0) node [below] {$\frac{3}{5}$};
\draw[blue] (\deb+8+0.15,0) node [below] {$\frac{4}{5}$};
\draw[blue] (\deb+2+0.15,-3) node [above] {$\frac{1}{5}$};
\draw[blue] (\deb+4+0.15,-3) node [above] {$\frac{2}{5}$};
\draw[blue] (\deb+6+0.15,-3) node [above] {$\frac{3}{5}$};
\draw[blue] (\deb+8+0.15,-3) node [above] {$\frac{4}{5}$};
\end{footnotesize}
\draw[->,>=latex] (E) to[out=-90, in =100, distance=1cm] (E');
\draw[->,>=latex] (F) to[out=-90, in =80,  distance=1.25cm] (F');
\draw[->,>=latex] (G) to[out=-90, in =90] (G');

\draw[color=black] (10,-6) -- (12,-6);
\draw (10,-6.5) -- (10,-5.5);
\draw (12,-6.5) -- (12,-5.5);
\draw[shift={(11,-6)},color=black] (0pt,2pt) -- (0pt,-2pt);
\draw (10.05,-6.45)--(10,-6.4)--(10.05,-6.35);
\draw (10.95,-6.45)--(11,-6.4)--(10.95,-6.35);
\draw (10,-6.4)--(11,-6.4);
{\footnotesize
\draw (10.6,-7) node {$\alpha_5(S)$};
}
\draw (13,-6) node {$Y_5(S)$};

\draw[color=black] (\deb,-6) -- (\deb+2,-6);
\draw (\deb,-6.5) -- (\deb,-5.5);
\draw (\deb+2,-6.5) -- (\deb+2,-5.5);
\foreach \x in {1.8,1.2,0.8}
\draw[shift={(\deb+\x,-6)},color=black] (0pt,2pt) -- (0pt,-2pt);
\draw (\deb+1.85,-6.45)--(\deb+1.8,-6.4)--(\deb+1.85,-6.35);
\draw (\deb+1.95,-6.45)--(\deb+2,-6.4)--(\deb+1.95,-6.35);
\draw (\deb+1.8,-6.4)--(\deb+2,-6.4);
{\footnotesize
\draw (\deb+1.9,-7) node {$\alpha_5(R)$};
}
\draw (\deb+3,-6) node {$Y_5(R)$};

\end{tikzpicture}
\caption{Example where $\mathcal{U}_{\epsilon}^q$ is not a neighbourhood of $S$}\label{counterexUalpha}
\end{figure}
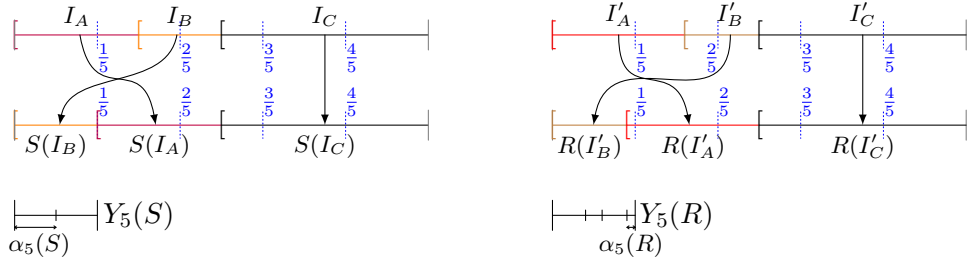

Then $\alpha_q(R)<\eta<\epsilon< \frac{1}{2q}=\alpha_q(S)$, so $R\notin\mathcal{U}_\eta^q$ though $R$ is at distance less than $2\eta$ from $S$.
\end{example}

\begin{remark} The condition (written in Remark \ref{U_ngb}) on $S$ corresponds to being outside a finite union of hyperplanes. This means that for every $S$ in a full-measure set $\mathcal{A}$, for every $\epsilon<\alpha_q(S)$, $\mathcal{U}_\epsilon^q(S)$ is a neighbourhood of $S$.
\end{remark}

\begin{lem}
Let $q\in\N$, $\epsilon>0$ and $S$ be an IET in $\mathcal{U}_{\epsilon}^q$.
Then:
$$ S\left(\mathcal{N}_{\epsilon}\left(\Delta(S) \cup \frac{1}{q}\N \right)\right) \subset \mathcal{N}_{\epsilon}(X_q(S)).$$
\end{lem}

\begin{proof}
Note that $  S(\mathcal{N}_{\epsilon}(\Delta(S)))=\mathcal{N}_{\epsilon}(S(\Delta(S))\cup \{1\})$ and that $S(\mathcal{N}_{\epsilon}(\frac{1}{q}\N))\subset \mathcal{N}_{\epsilon}(S(\frac{1}{q}\N))$.
Hence:
\begin{align*}
S\left(\mathcal{N}_{\epsilon}\left(\Delta(S) \cup \frac{1}{q}\N \right)\right) 
&\subset \mathcal{N}_{\epsilon}(S(\Delta(S))) \cup \mathcal{N}_{\epsilon}(\{1\})\cup \mathcal{N}_{\epsilon} \left( S\left(\frac{1}{q}\N\right)\right)  \\
&\subset \mathcal{N}_{\epsilon}(\underset{=\Delta(S^{-1})\cup\{0,1\}}{\underbrace{S(\Delta(S))\cup\{0,1\}}}) \cup \mathcal{N}_{\epsilon} \left( S\left(\frac{1}{q}\N\right)\right) \\
&\subset \mathcal{N}_{\epsilon}(X_q(S)).
\end{align*}
\end{proof}

\section{Proof of the main result}
We are now ready to prove Theorem \ref{thJ}. Let us start by building an IET with small support.

\subsection{Building an IET with small support}\label{section_small_supp}

The lemmas \ref{small_trans_length_bis} and \ref{lem_support} are very similar to the ones of \cite{DFG}, with small changes to overcome the non $q$-rationality of $S$.

\begin{lem}[5.5 of \cite{DFG}]\label{lemma5.6}
For all $\epsilon >0$, $q\in \N$, $m\in\N$ there exists $\eta >0$ such that if $S_0$ and $T_0$ are $q$-rational, if $w(s,t)$ is a word of length at most $m$ in the letters $s^{\pm 1}$ and $t^{\pm 1}$ such that $w(S_0,T_0)=\mathrm{id}$, and if $S$ and $T$ are $\eta$-close to $S_0$ and $T_0$ respectively, then
on each interval of $[0,1)\setminus \mathcal{N}_\epsilon(\frac{1}{q}\N)$,
the transformation $w(S,T)$ induces a translation of length strictly less than $\epsilon$.
\end{lem}

We use this lemma with $m=q!$ and $w=t^{q!}$ to prove the next lemma. We take $\eta < \frac{1}{2m}$  as in the proof of Lemma 5.5 in \cite{DFG}.
\begin{lem}\label{small_trans_length_bis}
For all $\epsilon >0$, $q\in \N$, there exists $\eta >0$ such that if $T_0$ is $q$-rational and if $T$ is $\eta$-close to $T_0$, then
for any IET $S$, the IET
$ST^{q!}S^{-1}$ induces translations on $[0;1) \, \setminus \, S(\mathcal{N}_{\epsilon}(\Delta(S) \cup \frac{1}{q}\N ))$ with translation lengths strictly less than $\epsilon$.
\end{lem}

\begin{proof}
Let $\eta < \frac{\epsilon}{2q!}$ and
 $x \notin S(\mathcal{N}_{\epsilon}(\Delta(S) \cup \frac{1}{q}\N ))$.

Then $S^{-1}(x) \notin \mathcal{N}_{\epsilon}(\Delta(S) \cup \frac{1}{q}\N)$. By Lemma \ref{lemma5.6}, $T^{q!}$ is a translation on each interval of $[0;1)  \, \setminus  \, \mathcal{N}_{\epsilon}(\frac{1}{q}\N)$ with translation length strictly less than $\epsilon$, thus $|T^{q!}S^{-1}(x)-S^{-1}(x)| < \epsilon$. Since $d(S^{-1}(x),\Delta(S)) \geqslant \epsilon$, this tells us that  $T^{q!}S^{-1}(x)$ is in the same interval of continuity of $S$ as $S^{-1}(x)$. Applying the local isometry $S$ to 
$| T^{q!}S^{-1}(x) - S^{-1}(x)| < \epsilon$ yields $|ST^{q!}S^{-1}(x) -x| < \epsilon$.
\end{proof}

\begin{lem}\label{lem_support}
For all $\epsilon >0$, $q\in \N$, set $\eta = \frac{\epsilon}{4q!}$. If $T_0$ is $q$-rational and if $T$ is $\eta$-close to $T_0$, then
for all IET $S$, the IET
$[T^{q!},ST^{q!}S^{-1}]$ induces the identity on each interval of $[0;1) \, \setminus \, \mathcal{N}_{\epsilon}(X_q(S))$.
\end{lem}

\begin{proof}
Let us write $A_{\epsilon}=[0;1) \, \setminus \, \mathcal{N}_{\epsilon}(X_q(S))$ and $A_{\frac{\epsilon}{2}}=[0;1) \, \setminus \, \mathcal{N}_{\frac{\epsilon}{2}}(X_q(S))$.

Note that $\mathcal{N}_{\frac{\epsilon}{2}}(\frac{1}{q}\N) \subset \mathcal{N}_{\epsilon}(X_q(S))$ and  $S(\mathcal{N}_{\frac{\epsilon}{2}}(\Delta(S)\cup\frac{1}{q}\N)) \subset \mathcal{N}_{\epsilon}(X_q(S))$.

Apply Lemma \ref{small_trans_length_bis} with $\frac{\epsilon}{2}$: whenever $T$  is closer than $\eta = \frac{\epsilon}{4q!}$ to $T_0$, the IET $ST^{q!}S^{-1}$ induces translations of translation lengths strictly smaller than  $\frac{\epsilon}{2}$ on each connected component of $[0;1) \, \setminus \, S(\mathcal{N}_{\frac{\epsilon}{2}}(\Delta(S)\cup\frac{1}{q}\N)) \supset A_{\frac{\epsilon}{2}}$. Moreover,
the IET $T^{q!}$ induces translations of translation lengths strictly smaller than  $\frac{\epsilon}{2}$ on each connected component of $[0;1) \, \setminus \, \mathcal{N}_{\frac{\epsilon}{2}}(\frac{1}{q}\N) \supset A_{\frac{\epsilon}{2}}$ (Lemma 5.5 of \cite{DFG}).

Let $x \in A_{\epsilon}$. Denote by $I=[a_1;a_2]$ the connected component of $x$ in $A_{\frac{\epsilon}{2}}$ and by $t, t' \in [-\frac{\epsilon}{2};\frac{\epsilon}{2}]$ the translation lengths on $I$ of $T^{q!}$ and $ST^{q!}S^{-1}$ respectively.
Both $x+t$ and $x+t'$ are in $I$ so $T^{q!}\circ ST^{q!}S^{-1} (x) = T^{q!} (x + t ') = x + t' + t$ and $ ST^{q!}S^{-1} \circ T^{q!}(x) = ST^{q!}S^{-1} (x + t ) = x + t + t' $.

Hence $T^{q!}$ and $ ST^{q!}S^{-1} $ commute on $A_{\epsilon} = [0;1) \, \setminus \, \mathcal{N}_{\epsilon}(X_q(S))$. We conclude that $[T^{q!},ST^{q!}S^{-1}]$ induces the identity on each interval of $[0;1) \, \setminus \, \mathcal{N}_{\epsilon}(X_q(S))$.
\end{proof}

In other words, the IET $U=[T^{q!},ST^{q!}S^{-1}]$ has its support included in $\mathcal{N}_{\epsilon}(X_q(S))$ and thus also in $\mathcal{N}_{\epsilon}(Z_q(S))$. 

This means that the support of $U$ is both small and included in the "periodic" set $\mathcal{N}_{\epsilon}(Z_q(S))= \mathcal{N}_{\epsilon}(Y_q(S)) \sqcup \left(\mathcal{N}_{\epsilon}(Y_q(S)) + \frac{1}{q}\right) \sqcup ... \sqcup \left(\mathcal{N}_{\epsilon}(Y_q(S)) +\frac{q-1}{q} \right)$. This periodicity is important since the only control we have on the drift (as we will see in section \ref{section_drift_supp}) is modulo $\frac{1}{q}$.

\subsection{Drifting the support} \label{section_drift_supp}

\begin{definition}[drifting vector and drifting direction] Let $\sigma\in\mathfrak{S}_n$, and $T\in \IET_\sigma$. Write $(l_1,...,l_n)=\lambda(T)$
and $(t_1,...,t_n)=\omega(T)$.

We say that $\sigma$ is \emph{driftable} if there exist $(d_1,...,d_n)$, where $\sum d_i = 0$, such that the IET $T'$ defined by $\lambda(T')=(l_1+d_1,...,l_n+d_n)$ and $\pi(T') =\pi(T)$ has translation lengths $\omega(T')=(t_1+v_1,...,t_n+v_n)$ with all $v_i$ positive.

In this case, $(d_1,...,d_n)$ is called a \emph{drifting direction} and $(v_1,...,v_n)$ a \emph{drifting vector}.
\end{definition}

\begin{remark} A permutation is driftable if and only if it is admissible. See Proposition 5.12 of \cite{DFG} for a proof. We highlight the fact that driftability depends only on the permutation $\sigma$ of the IET and neither on its translation lengths nor on its length of intervals of continuity. 
\end{remark}

\begin{definition}
Let $T_0$ be an IET and $\sigma=\pi(T_0)\in\mathfrak{S}_n$ be its underlying permutation.
Let $u\in\R^n$.
We define $T_{\theta}^u \in \IET_{\sigma}$ by $\lambda(T_{\theta}^u) = \lambda(T_0) + \theta u$, where $\theta$ is small enough to ensure that all the lengths $l_i(T_{\theta}^u)$ of the intervals of continuity of $T_\theta^u$ are positive.
\end{definition}

\begin{remark}
One can define $T_\theta^u$ even if $\sigma = \pi(T)$ is not driftable.
\end{remark}

\begin{definition}
Let $u = (u_1,...,u_n)\in\R^n$. We note $u_\text{min} = \min\{u_1,...,u_n\}$ and $u_\text{max} = \max\{u_1,...,u_n\}$.
\end{definition}

\begin{lem}\label{power_bis}
Let $\sigma \in \mathfrak{S}_n$ be a driftable permutation, let $d \in \R^n$ be a drifting direction and $v\in\R^n$ be the associated drifting vector.
Let $q \in \N$ and $T_0 \in \IET_{\sigma}([0;1))$ be $q$-rational.
Let $\rho = \frac{v_{\text{max}}}{v_{\text{min}}}$, $\alpha<\frac{1}{2q}$, $\epsilon < \frac{\alpha}{11\rho}$ and $\mu < \frac{\theta v_{\text{min}}}{4}$.
Let  $\theta < \frac{\epsilon}{v_{\text{min}}}$ and small enough so that $T_{\theta}^{d}$ is well defined. 
Let $T$ be $\mu$-close to $T_{\theta}^{d}$. 

Then there exists $k \in \N$ such that all the translation lengths of $T^k$ are in $[2\epsilon; \alpha-2\epsilon ] \text{ mod } \frac{1}{q}$.
\end{lem}

\begin{proof}
First, we have 
$$\theta v_{\text{min}} < \epsilon < \frac{\alpha}{11\rho} \leqslant \frac{1}{22q\rho} < \frac{1}{q}$$
and
$$
2\mu < \frac{1}{2}\theta v_{\text{min}} < \theta v_{\text{min}}
$$
and
$$
\theta v_{\text{max}} + 2\mu < \theta v_{\text{max}} + \theta v_{\text{min}} \leqslant 2 \theta v_{\text{max}} = 2 \theta v_{\text{min}} \rho < 2\epsilon \rho < \frac{1}{11q} <\frac{1}{q} .
$$

We deduce (Lemma 5.9 of \cite{DFG}) that all the translation lengths of $T$ are in the interval $[ \theta v_{\text{min}} -2\mu ; \theta v_{\text{max}} + 2\mu ] \text{ mod } \frac{1}{q}$, which is included in $ [ \frac{1}{2}\theta v_{\text{min}} ; \frac{3}{2} \theta v_{\text{max}} ] $.

For every $k \in\N$ such that $\frac{3k}{2} \theta v_{\text{max}} < \frac{1}{q}$, the IET $T^k$ has all of its translation lengths in $ [ \frac{k}{2}\theta v_{\text{min}} ; \frac{3k}{2} \theta v_{\text{max}} ] \text{ mod } \frac{1}{q} $.
We want an integer $k$ such that both $ 2\epsilon < \frac{k}{2}\theta v_{\text{min}} $ and $ \frac{3k}{2} \theta v_{\text{max}} < \alpha - 2\epsilon  \,$.
Take $k$ such that $ 4\epsilon < k \theta v_{\text{min}} \leqslant 5 \epsilon$ (this is possible because $\theta v_{\text{min}} < \epsilon$).
Then we have:
$$
 \frac{3k}{2} \theta v_{\text{max}} =  \frac{3k}{2} \theta v_{\text{min}} \rho \leqslant \frac{3}{2}\times 5 \epsilon \rho < \frac{15}{2} \frac{\alpha}{11} = \alpha - \frac{7\alpha}{22} < \alpha - \frac{2\alpha}{11} < \alpha - 2\epsilon.
$$
Thus $T^k$ has all its translation lengths in $ [2\epsilon ; \alpha-2\epsilon]  \text{ mod } \frac{1}{q}$.
\end{proof}

\subsection{Conclusion} \label{ccl} We are 
ready to build a nontrivial relation.

\begin{prop}\label{not_free_bis}Assume $T_0$ to be $q$-rational and its underlying permutation $\sigma \in \mathfrak{S}_n$ to be driftable. Let $S$ be any IET.

Then there exists a set $\mathcal{U} \subset \IET([0;1))$ that contains $S$ and an open set $\mathcal{V} \subset \IET_{\sigma}([0;1))$ that accumulates on $T_0$ such that $\langle R,T \rangle$ is not free of rank 2 whenever $R \in \mathcal{U}$ and $ T \in \mathcal{V}$. 

Moreover if $S$ is rational or such that for every distinct $x,y\in\Delta(S)$, one has $x +\frac{1}{q}\N \neq y+\frac{1}{q}\N$, then the set $\mathcal{U}$ is a neighbourhood of $S$. 
\end{prop}

\begin{proof}
If $S$ is rational, then Proposition 8 of \cite{DFG} proves the previous statement. Therefore we now assume that $S$ is not rational.

Let $\delta < \alpha_q(S)$. Note that $\alpha_q(S) \leqslant \frac{1}{2q}$ because $S$ is not $q$-rational. 

Let $d$ be a drifting direction for $\sigma$ and $v$ be the associated drifting vector.
Let $\rho = \frac{v_{\text{max}}}{v_{\text{min}}}$ and take

\begin{center}
\begin{tabular}{cccc}
$\epsilon < \frac{\delta}{11\rho}$; &
$\eta < \frac{\epsilon}{4q!}$;&
$\theta < \text{min}(\frac{\epsilon}{v_{\text{min}}} , \frac{\eta}{2 ||d||_1})$;&
$\mu < \text{min}(\frac{\theta v_{\text{min}}}{4}, \frac{\eta}{2})$.
\end{tabular}
\end{center}

Take $\mathcal{U} = \mathcal{U}_{\delta}^q$. If $S$ is such that for every distinct $x,y\in\Delta(S)$, one has $x +\frac{1}{q}\N \neq y+\frac{1}{q}\N$, then the set $\mathcal{U}$ is a neighbourhood of $S$ (Remark \ref{U_ngb}).

Take $T_{\theta}^d \in \IET_{\sigma}$ such that $\lambda(T_{\theta}^d) = \lambda(T_0) + \theta d$.
Let $\mathcal{V}_{d, \theta}$ be the set of IETs that are  at distance strictly less than $\mu$ from $T_{\theta}^d$. 
Then every $T \in \mathcal{V}_{d, \theta}$ satisfies:
$$
d(T,T_0) \leqslant d(T,T_{\theta}^d) +	d(T_{\theta}^d,T_0) < \mu + \theta ||d||_1 < \frac{\eta}{2} + \frac{\eta}{2} = \eta.
$$

From Lemma \ref{lem_support} 
we know that for every $(R,T)\in\mathcal{U}\times\mathcal{V}_{d, \theta}$ the IET $U = [T^{q!},RT^{q!}R^{-1}] $ has its support included in $\mathcal{N}_{\epsilon}(Z_q(R))$.
By the definition of $\alpha_q(R)$, the smallest connected component of $[0;1) \ \setminus \  \mathcal{N}_{\epsilon}(Z_q(R))$ has length at least $\alpha_q(R) - 2\epsilon$, hence at least $ \delta - 2\epsilon$ because $R\in\mathcal{U}=\mathcal{U}_{\delta}^q$.
It is therefore enough to drift the support of $U$ with a drift in $[2\epsilon , \delta - 2\epsilon]$ modulo $\frac{1}{q}\Z$. From Lemma \ref{power_bis}, there exists $k\in\N$ such that $[U,T^{k}UT^{-k}] = \mathrm{id}$.

It remains to check that this relation is not trivial. Denote by $u=t^{q!} r t^{q!} r^{-1}t^{-q!} r t^{-q!} r^{-1}$ the (non-trivial) word over the letters $r^{\pm 1}$ and $t^{\pm 1}$ such that $u(S,T) =U$. 
The word $w=ut^{k}ut^{-k}u^{-1}t^{k}u^{-1}t^{-k}$ is equal to:
\begin{align*}
(t^{q!} r t^{q!} r^{-1}t^{-q!} r t^{-q!} r^{-1})
.t^{k}.
(t^{q!} r t^{q!} r^{-1}t^{-q!} r t^{-q!} r^{-1})
.t^{-k}. 
 \qquad  \qquad \qquad
\\
 \qquad  \qquad  \qquad
(r t^{q!} r^{-1} t^{q!} r t^{-q!} r^{-1} \textcolor{red}{t^{-q!}})
.\textcolor{red}{t^{k}}.
(r t^{q!} r^{-1} t^{q!} r t^{-q!} r^{-1} t^{-q!})
.t^{-k} 
\end{align*}
which is reduced, except if $k=q!$ in which case it reduces to the reduced word:
$$
t^{q!} r t^{q!} r^{-1}t^{-q!} r t^{-q!} r^{-1}
t^{k}
t^{q!} r t^{q!} r^{-1}t^{-q!} r t^{-q!} r^{-1}
t^{-k}
r t^{q!} r^{-1} \textcolor{red}{t^{2q!}} r t^{-q!} r^{-1} t^{-q!}
t^{-k}. 
$$
So the relation is not trivial.
Hence for every $(R,T)\in\mathcal{U}\times\mathcal{V}_{d, \theta}$ the subgroup $\langle R,T \rangle$ is not free of rank 2.

In order to have an open set $\mathcal{V}$ that accumulates on $T_0$ we take the union of all the convenient $\mathcal{V}_{d,\theta}$ (which also depend on $\delta$). More precisely, we define:
$$\mathcal{V}_{d} =
\underset{\epsilon < \frac{\delta}{11\rho}}{\bigcup} \, \,
\underset{\eta < \frac{\epsilon}{4q!}}{\bigcup} \, \,
\underset{\theta < \text{min}(\frac{\epsilon}{v_{\text{min}}} , \frac{\eta}{2 ||d||_1})}{\bigcup} \mathcal{V}_{d, \theta}
\, \text{ and } \,
\mathcal{V} = \underset{d \text{ drifting direction}}{\bigcup} \mathcal{V}_{d}.
$$
\end{proof}

\begin{proof}[Proof of Theorem \ref{thJ}]
Let $S$ be an IET on $[0;1)$.
For every $q$-rational driftable IET $T_0$, fix $\delta < \alpha_q(S)$ and denote by $\mathcal{V}(T_0)$ the open set given by the previous proposition and define:
$$
\Omega_a(S) =
\underset{q\geqslant 2}{\bigcup} \quad
\underset{2\leqslant n \leqslant q}{\bigcup} \quad
\underset{\sigma \in \mathfrak{S}_n \, \text{driftable}}{\bigcup} \quad
\underset{T_0 \in \IET_{\sigma} \, q-\text{rational}}{\bigcup}  \mathcal{V}(T_0).
$$

The set $\Omega_a(S)$ is a union of open sets and hence is open, it is dense in $ \IET_a([0;1))$ because of the density of $\Q$ in $\R$. And the proof of the previous proposition shows that for every $T \in \Omega_a(S)$, $\langle S,T \rangle$ is not free of rank 2.

\end{proof}

\subsection{Looking for smaller relations}

We have shown that given any IET $S$, there is a dense subset of admissible IETs sharing a relation with $S$. The next natural question is the measure of this set (in each connected component of $\IET([0;1))$). It is non-zero but seems very small. Indeed
the open set $\mathcal{V}(T_0)$ that we built around a $q$-rational IET $T_0$ has a diameter which is bounded by $\frac{1}{q!q}$ times a constant (depending only on the underlying permutation $\sigma$ of $T_0$).
For instance, there does not exist a $T_0$ such that the Arnoux-Yoccoz example is in $\mathcal{V}(T_0)$. This example is interesting because it is self-similar (applying Rauzy induction a certain number of times gives the same IET).

But maybe we could have built a bigger set $\Omega(S)$. We briefly discuss this issue in the following subsection.

\subsubsection{Arnoux-Yoccoz's example} \label{AYex}

Let us present an example by Arnoux and Yoccoz (see \cite{AY} and \cite{Arnoux}).
Let $a$ be the only real number such that $a^3+a^2+a=1$.
Define the IET $g$ by 
$$
\lambda(g)=\left(\lambda_A=\frac{a}{2}, \lambda_{A'}=\frac{a}{2}, \lambda_B=\frac{a^2}{2},\lambda_{B'}=\frac{a^2}{2}, \lambda_{C}=\frac{a^3}{2}, \lambda_{C'}=\frac{a^3}{2}\right)
$$
and 
$$
\pi(g)= \left(
\begin{array}{cccccc}
A&A'&B&B'&C&C'\\
A'&A&B'&B&C'&C
\end{array}
\right)
.$$
Let $h$ be the rotation of angle $\frac{1}{2}$, namely $\lambda(h)= (\lambda_A = \frac{1}{2},\lambda_B=\frac{1}{2})$ and $\pi(h)=\left(
\begin{array}{cc}
A&B\\
B&A
\end{array}
\right)$.
Finally, define $f=h \circ g$. We have:
\begin{small}
$$
\lambda(f)=\left(\lambda_A=\frac{1-a}{2}, \lambda_B=a-\frac{1}{2}, \lambda_C=\frac{a}{2},\lambda_D=\frac{a^2}{2}, \lambda_E=\frac{a^2}{2}, \lambda_F=\frac{a^3}{2},\lambda_G=\frac{a^3}{2}\right)
$$
\end{small}
and 
$$
\pi(f)= \left(
\begin{array}{ccccccc}
A&B&C&D&E&F&G\\
B&E&D&G&F&C&A
\end{array}
\right)
.$$

\begin{figure}[ht]
\center
\begin{tikzpicture}[x=0.105cm,y=0.4cm]
\def \a{27}
\def \bb{4}
\def \b{15}
\def \c{8}
\def \L{2*\a+2*\b+2*\c}

\def \cola{blue}
\def \colb{Turquoise}
\def \colc{red}
\def \cold{Maroon}
\def \cole{Orange}
\def \colf{violet}
\def \colg{magenta}

\foreach \h in {0,-3}
{\draw[gray] (0,\h+-0.5) -- (0,\h+0.5);
\draw[gray] (\L,\h+-0.5) -- (\L,\h+0.5);}
\draw [\cola] (0,0) -- (\a-\bb,0);
\draw [\cola,double] (\a-\bb,0) -- (\a,0);
\draw [\colc] (\a-\bb,0) -- (\a,0);
\draw [\colb] (\a,0) -- (2*\a,0);
\draw [\cold] (2*\a,0) -- (2*\a+\b,0);
\draw [\cole] (2*\a+\b,0) -- (2*\a+2*\b,0);
\draw [\colf] (2*\a+2*\b,0) -- (2*\a+2*\b+\c,0);
\draw [\colg] (2*\a+2*\b+\c,0) -- (2*\a+2*\b+2*\c,0);
\draw [\colb] (0,-3) -- (\a,-3);
\draw [\cola] (\a,-3) -- (2*\a-\bb,-3);
\draw [\cola,double] (2*\a-\bb,-3) -- (2*\a,-3);
\draw [\colc] (2*\a-\bb,-3) -- (2*\a,-3);
\draw [\cole] (2*\a,-3) -- (2*\a+\b,-3);
\draw [\cold] (2*\a+\b,-3) -- (2*\a+2*\b,-3);
\draw [\colg] (2*\a+2*\b,-3) -- (2*\a+2*\b+\c,-3);
\draw [\colf] (2*\a+2*\b+\c,-3) -- (2*\a+2*\b+2*\c,-3);
\draw (50,-3) node {$[$} ;
\foreach \x in {0,\a,2*\a,2*\a+\b,2*\a+2*\b,2*\a+2*\b+\c}
\draw[shift={(\x,0)}] node {$[$} ;
\draw[gray] (\a-\bb,0) node {$[$} ;
\foreach \x in {\a,2*\a,2*\a+\b,2*\a+2*\b,2*\a+2*\b+\c}
\draw[shift={(\x,-3)},gray] node {$[$} ;
\draw (0,-3) node {$[$} ;

{\footnotesize
\draw (\a/2,0) node (A) {};
\draw (\a/2-\bb/2,0) node [above]{$I_A$};
\draw (\a*3/2,0) node (B) {};
\draw (\a-\bb/2,0) node [above]{$I_B$};
\draw (\a*3/2,0) node [above] {$I_C$};
\draw (2*\a+\b/2,0) node (C) {};
\draw (2*\a+\b/2,0) node [above]{$I_D$};
\draw (2*\a+\b*3/2,0) node (D) {};
\draw (2*\a+\b*3/2,0) node [above]{$I_E$};
\draw (2*\a+2*\b+\c/2,0) node (E) {};
\draw (2*\a+2*\b+\c/2,0) node [above]{$I_F$};
\draw (2*\a+2*\b+\c*3/2,0) node (F) {};
\draw (2*\a+2*\b+\c*3/2,0) node [above]{$I_G$};

\draw (\a*3/2,-3) node (A'){} ;
\draw (\a*3/2,-3) node [below]{$g(I_A)$};
\draw (\a/2,-3) node (B') {};
\draw (\a/2,-3) node [below]{$g(I_C)$};
\draw (2*\a-\bb/2,-3) node [below]{$g(I_B)$};
\draw (2*\a+\b*3/2,-3) node (C') {};
\draw (2*\a+\b*3/2,-3) node [below]{$g(I_D)$};
\draw (2*\a+\b/2,-3) node (D') {};
\draw (2*\a+\b/2,-3) node [below]{$g(I_E)$};
\draw (2*\a+2*\b+\c*3/2,-3) node (E') {};
\draw (2*\a+2*\b+\c*3/2+0.15,-3) node [below]{$g(I_F)$};
\draw (2*\a+2*\b+\c/2,-3) node (F') {};
\draw (2*\a+2*\b+\c/2-0.15,-3) node [below]{$g(I_G)$};
}

\draw[->,>=latex] (A) to[out=-90, in =90,distance=0.8cm] (A');
\draw[->,>=latex] (B) to[out=-90, in =90,distance=0.8cm] (B');
\draw[->,>=latex] (C) to[out=-90, in =90] (C');
\draw[->,>=latex] (D) to[out=-90, in =90] (D');
\draw[->,>=latex] (E) to[out=-90, in =90] (E');
\draw[->,>=latex] (F) to[out=-90, in =90] (F');

\draw [->] (-0.15*5,0) to [out = -135,in =135] (-0.15*5,-3);
\draw (-0.9*5,-3/2) node {$g$}; 

\draw [->,dashed] (-0.15*5,-3) to [out = -135,in =135] (-0.15*5,-6);
\draw (-0.9*5,-3-3/2) node {$h$};

\def \a{27}
\def \b{23}
\def \c{4}
\def \d{15}
\def \e{8}
\def \L{\a+\b+\c+2*\d+2*\e}
\def \h{6}

\draw[gray] (0,-\h-0.5) -- (0,-\h+0.5);
\draw[gray] (\L,-\h-0.5) -- (\L,-\h+0.5);
\draw (0,-3) -- (\L,-3);
\draw [\cola,double] (0,-\h) -- (\c,-\h);
\draw [\colc] (0,-\h) -- (\c,-\h);
\draw [\cole] (\c,-\h) -- (\c+\d,-\h);
\draw [\cold] (\c+\d,-\h) -- (\c+2*\d,-\h);
\draw [\colg] (\c+2*\d,-\h) -- (\c+2*\d+\e,-\h);
\draw [\colf] (\c+2*\d+\e,-\h) -- (\c+2*\d+2*\e,-\h);
\draw [\colb] (\c+2*\d+2*\e,-\h) -- (\c+2*\d+2*\e+\a,-\h);
\draw [\cola] (\c+2*\d+2*\e+\a,-\h) -- (\c+2*\d+2*\e+\a+\b,-\h);
\foreach \x in {0,\c,\c+\d,\c+2*\d,\c+2*\d+\e,\c+2*\d+2*\e,\c+2*\d+2*\e+\a}
\draw[shift={(\x,-\h)}] node {$[$} ;
\draw (50,-\h) node {$[$} ;

\draw[->,>=latex,dashed] (25,-3-0.3) to[out=-90, in =90, distance = 0.9cm] (75,-\h+0.3);
\draw[->,>=latex,dashed] (75,-3-0.3) to[out=-90, in =90, distance = 0.9cm] (25,-\h+0.3);
{\footnotesize
\draw (\c/2,-\h) node [below=3]{$f(I_B)$};
\draw (\c+\d/2,-\h) node [below=3]{$f(I_E)$};
\draw (\c+3*\d/2,-\h) node [below=3]{$f(I_D)$};
\draw (\c+2*\d+\e/2,-\h) node [below=3]{$f(I_G)$};
\draw (\c+2*\d+3*\e/2,-\h) node [below=3]{$f(I_F)$};
\draw (\L-\b-\a/2,-\h) node [below=3]{$f(I_C)$};
\draw (\L-\b/2,-\h) node [below=3]{$f(I_A)$};
}
\draw [->,double] (-0.15*5,0) to [out = -180,in =180] (-0.15*5,-6);
\draw (-1.9*5,-3) node {$f$}; 
\draw (\L+2.3*5,-3) node {}; 
\end{tikzpicture}
\caption{An example of the IET $f$ by Arnoux and Yoccoz}
\end{figure}
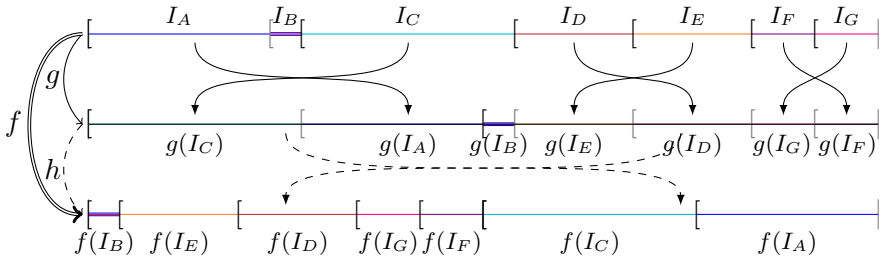

For every $X \in \{A,B,...,G\}$, $\lambda_X$ is an algebraic number of degree 3, so there exist $A(\lambda_X)$ such that for every rational number $\frac{p}{q}$, $\left| \lambda_X - \frac{p}{q}\right| \geqslant \frac{A(\lambda_X)}{q^3}$ (by Liouville's Theorem).
Moreover we easily compute such $A(\lambda_X)$. We conclude that a $q$-rational IET $T_0$ is far from $f$ by at least:
$$
d(f,T_0) \geqslant \sum_{X\in\{A,B,...,G\}} \frac{A(\lambda_X)}{q^3} \geqslant \frac{282}{q^2} \ > \ \frac{1}{q!}.
$$

On the other hand $f$ should be $\eta$-close to $T_0$ to be in $\mathcal{V}(T_0)$, where $\eta$ is less than $\frac{1}{q!}$. This proves that $f$ is in none of the $\mathcal{V}(T_0)$.

The factor $\frac{1}{q!}$ comes from the length of the basic relation we used ($T_0^{q!}=\mathrm{id}$).
To increase the size of $\mathcal{V}(T_0)$, we can take a smaller basic relation 
like $ST_0^{k}S^{-1}=\mathrm{id}$ where $k$ is the order of $T_0$.
This gives bigger open sets $\mathcal{V}(T_0)$ around $q$-rational IETs $T_0$ that have "small" order.

We want to know if it is enough. For every $q$ between 20 and 20000, we have computed the closest $q$-rational IET to $f$, its distance $\delta$ to $f$, its order \textcursive{o} and the bound  $\mathfrak{b} = 40q($\textcursive{o}$+2) \delta$.
The condition $\mathfrak{b} < 1$ is necessary so that $f$ is in an enlarged $\mathcal{V}(T_0)$.

\begin{figure}[h!]
\centering
\includegraphics[scale=0.2,trim=1cm 0cm 1cm 2.8cm, clip=true]{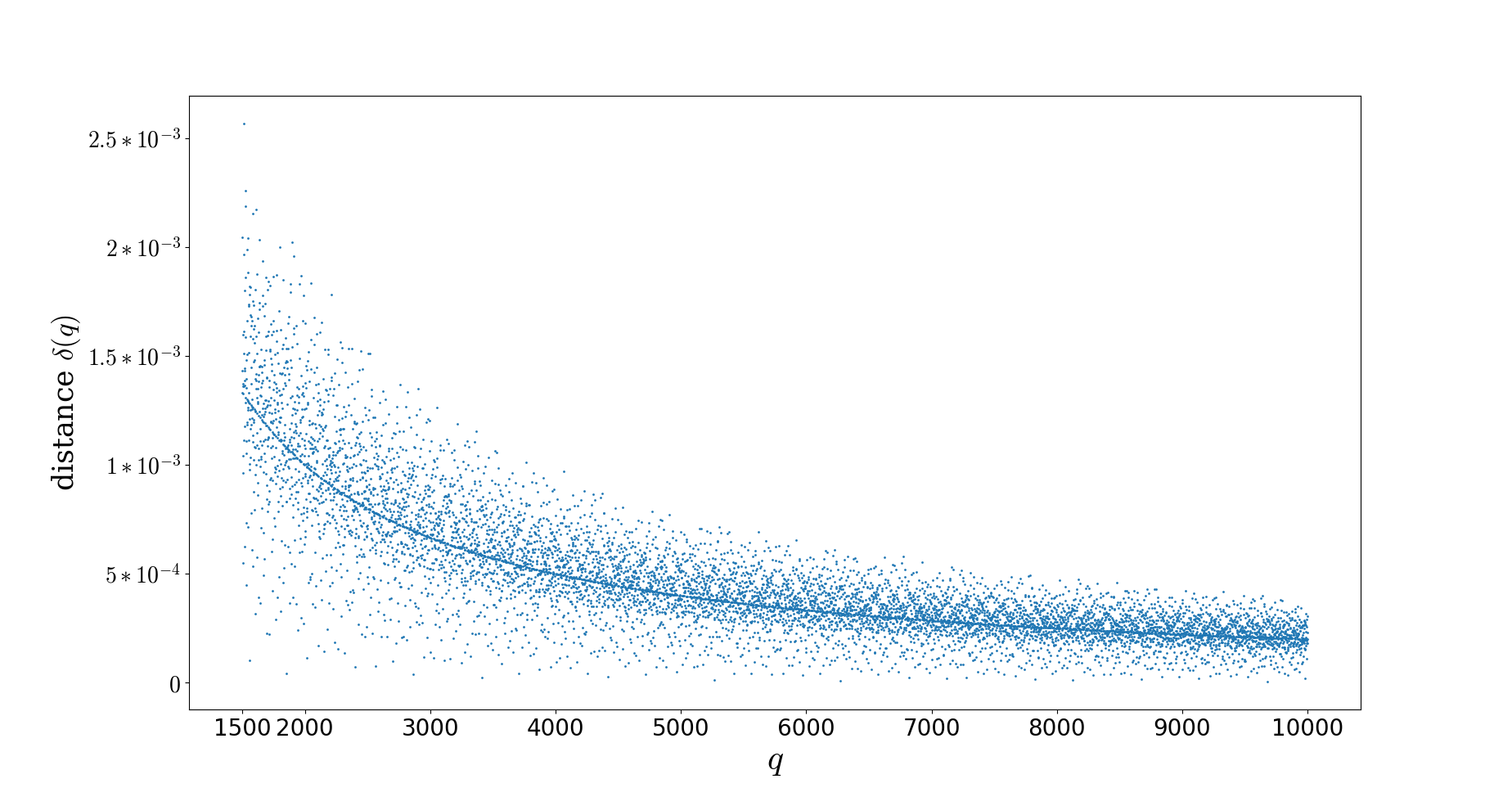}
\caption{\centering Distance $\delta(q)$ between $f$ and the closest $q$-rational IET to $f$}\label{figDist}
\end{figure}

The distance $\delta(q)$ between $f$ and the closest $q$-rational IET seems to decrease approximately like $\frac{1}{q}$ (and not faster). See Figure \ref{figDist}.

\begin{remark}
The distance between $f$ and the closest $q$-rational IET to $f$ is at most equal to $\frac{5}{q}$. Computations suggest that this distance does not decrease faster.
\end{remark}

\begin{figure}[h!]
\centering
\includegraphics[scale=0.2,trim=1cm 0cm 1cm 2.8cm, clip]{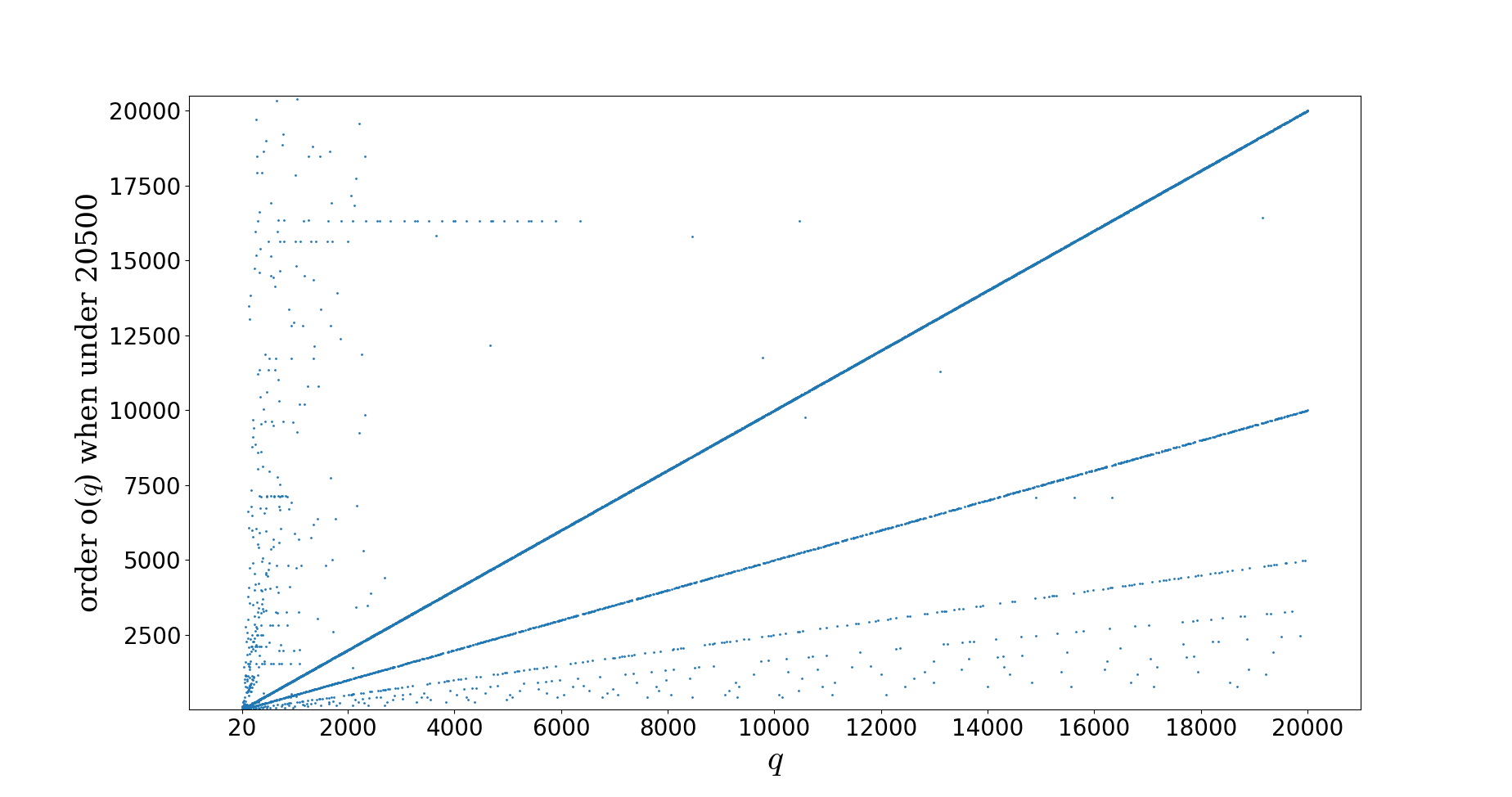}
\caption{Order \textcursive{o}$(q)$ of the closest $q$-rational IET to $f$}\label{figOrder}
\end{figure}

For some values of $q$, the order \textcursive{o}($q$) of this $q$-rational IET (the closest to $f$) increases linearly. See Figure \ref{figOrder}.

\begin{figure}[h!]
\centering
\includegraphics[scale=0.2,trim=1cm 0cm 1cm 2cm, clip]{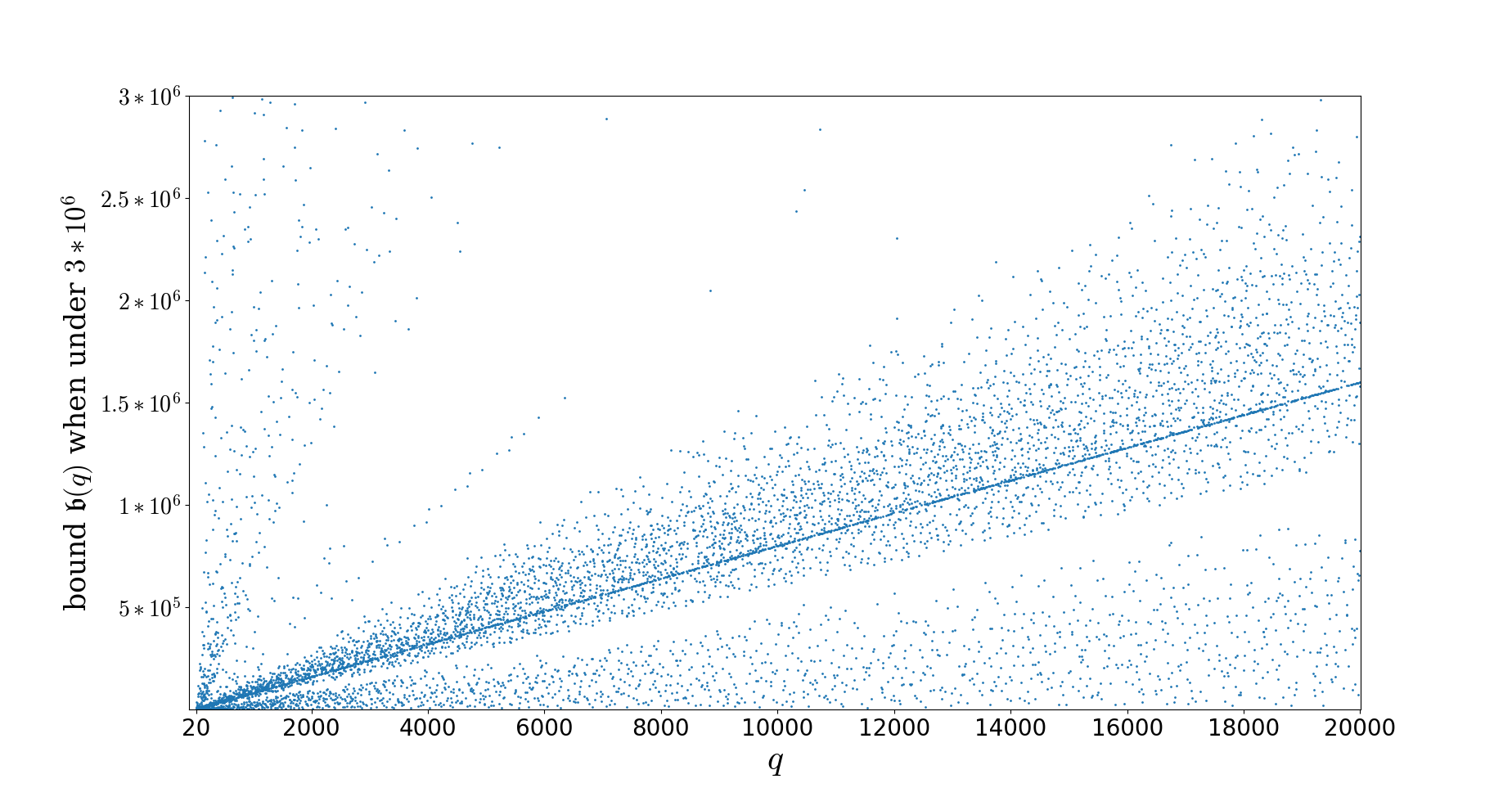}
\includegraphics[scale=0.2,trim=1cm 0cm 1cm 2cm, clip]{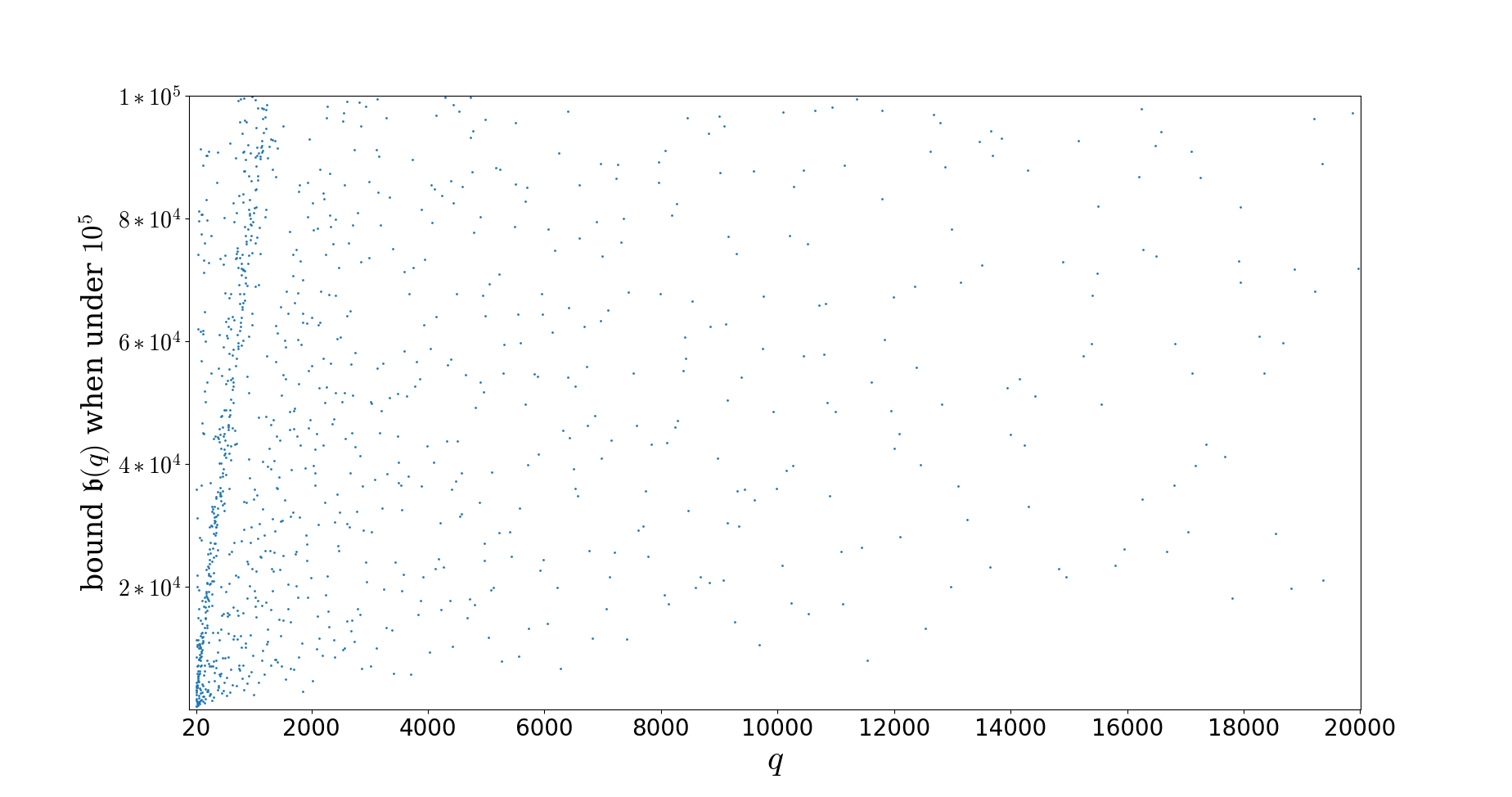}
\caption{Bound $\mathfrak{b}(q)$}\label{figBound}
\end{figure}

As a result, the bound $\mathfrak{b}(q)$ increases at least linearly and it seems that there is no rational IET $T_0$ close enough to $f$ so that   $f\in\mathcal{V}(T_0)$. See Figure \ref{figBound}.

\section{Free group in the group of affine interval exchange transformations}

\begin{definition}[AIET] \label{def_AIET} An \emph{affine interval exchange transformation} (AIET) on $[0;1)$ is a bijection from $[0;1)$ onto itself which is everywhere continuous on the right, continuous except on a finite number of points and differentiable except on a finite set and with constant differential on every interval where it is defined.

We denote by $\Delta(T)$ the union of discontinuities of $T$ and of the finite set where its differential $T'$ is not defined.
\end{definition}

In the following figures, we represent in dark (resp. light) colors the intervals that are expanded (resp. contracted) by the AIET and their images in light (resp. dark) colors.

Contrary to what happens in the group of IETs, we can easily find a free subgroup in the group of AIETs because its elements can contract and expand some parts of the domain of definition.

\begin{proof}[Proof of Proposition \ref{free_AIET}]
Let $V,W,X,Y$ be four disjoint subintervals of $[0;1)$ that do not cover $[0;1)$. Let $f$ be an AIET that sends $V$ onto $[0;1) \setminus W$ (and $[0;1)\setminus V$ onto $W$). This implies that $f(X)$ and $f(Y)$ are included in $W$ and $f^{-1}(X)$ and $f^{-1}(Y)$ are included in $V$. Let $g$ be an AIET that sends $X$ onto $[0;1) \setminus Y$ (and $[0;1)\setminus X$ onto $Y$). This implies that $g(V)$ and $g(W)$ are included in $Y$ and $g^{-1}(V)$ and $g^{-1}(W)$ are included in $X$.
See Figure \ref{fig_AIET}.
Applying the ping-pong lemma (See II. B. 24. of \cite{DeLaHarpe}), we see that $<f,g>$ is a free group of rank 2.
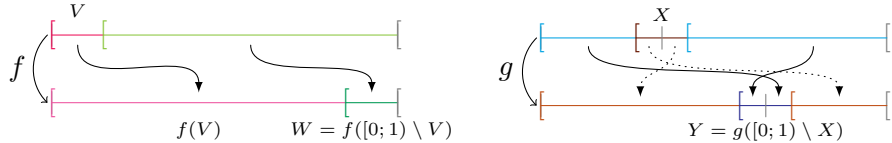
\begin{figure}[h!]
\center
\begin{tikzpicture}[line cap=round,line join=round,x=2.3cm,y=0.3cm]
\def \a{0.3}
\def \b{0.3}
\def \L{2}
\def \h{3}

\def \colabis{Rhodamine}
\def \cola{OrangeRed}
\def \colbbis{ForestGreen}
\def \colb{LimeGreen}

\draw [gray] (0,-0.5) -- (0,0.5);
\draw [gray] (\L,-0.5) -- (\L,0.5);
\draw[\cola] (0,0) -- (\a,0);
\draw[\cola] (0,0) node{$[$};
\draw[\colb] (\a,0) -- (\L,0);
\draw[\colb] (\a,0) node{$[$};
\draw [gray] (\L,0) node{$[$};

\draw [gray] (0,-\h-0.5) -- (0,-\h+0.5);
\draw [gray] (\L,-\h-0.5) -- (\L,-\h+0.5);
\draw[\colabis] (0,-\h) node{$[$};
\draw[\colabis] (0,-\h) -- (\L-\b,-\h);
\draw[\colbbis] (\L-\b,-\h) node{$[$};
\draw[\colbbis] (\L-\b,-\h) -- (\L,-\h);
\draw [gray] (\L,-\h) node{$[$};

\draw (\a/2,0) node (A)  {} ;
\draw (\a/2,0) [above = 3] node {\tiny $V$};
\draw (\L/2+\a/2,0) node (B) { };
\draw (\L/2-\b/2,-\h) node (A') { };
\draw (\L/2-\b/2,-\h)[below=3] node{\tiny $f(V)$};
\draw (\L-\b/2,-\h) node (B') { };
\draw (\L-\b/2,-\h) [below = 3] node{\tiny $W = f([0;1)\setminus V)$};
\draw[->,>=latex] (A) to[out=-90, in =90] (A');
\draw[->,>=latex] (B) to[out=-90, in =90,distance =15] (B');

\draw [->] (-0.05/2,0) to [out = -135,in =135] (-0.05/2,-\h);
\draw (-0.2,-\h/2) node {$f$}; 
\draw (\L+0.2,-\h/2) node {}; 
\end{tikzpicture}
\hspace{0.1cm}
\begin{tikzpicture}[line cap=round,line join=round,x=2.3cm,y=0.3cm]
\def \a{0.3}
\def \b{0.3}
\def \C{0.55} 
\def \c{0.3}
\def \D{1.15} 
\def \d{0.3}
\def \L{2}
\def \h{3}

\def \cola{Blue}
\def \colabis{Cerulean}
\def \colc{Sepia}
\def \colcbis{Bittersweet}

\draw [gray] (0,-0.5) -- (0,0.5);
\draw [gray] (\L,0) node{$[$};
\draw[\colabis] (0,0) -- (\C,0);
\draw[\colabis] (0,0) node{$[$};
\draw [\colc] (\C,0) node{$[$};
\draw [\colc] (\C,0) -- (\C+\c,0);
\draw [gray] (\C+\c/2,-0.5) -- (\C+\c/2,0.5);
\draw [\colabis] (\C+\c,0) -- (\L,0);
\draw [\colabis] (\C+\c,0) node{$[$};

\draw [gray] (0,-\h-0.5) -- (0,-\h+0.5);
\draw [gray] (\L,-\h) node{$[$};
\draw [\colcbis] (0,-\h) -- (\D,-\h);
\draw [\colcbis] (0,-\h) node{$[$};
\draw [\cola] (\D,-\h) node{$[$};
\draw[\cola] (\D,-\h) -- (\D+\d,-\h);
\draw [gray] (\D+\d/2,-\h-0.5) -- (\D+\d/2,-\h+0.5);
\draw [\colcbis] (\D+\d,-\h) -- (\L,-\h);
\draw [\colcbis] (\D+\d,-\h) node{$[$};

\begin{footnotesize}
\draw (\C/2,0) node (A) { };
\draw (\C+\c/4,0) node (B1) { };
\draw (\C+\c*3/4,0) node (B2) { };
\draw (\C+\c/2,0) node [above=3]{\tiny $X$};
\draw (\D/2+\L/2,0) node (C) { };
\draw (\D/2,-\h) node (B2') { } ;
\draw (\D+\d/4,-\h) node (A2') { };
\draw (\D+\d/2,-\h) node [below = 3]{\tiny $Y = g([0;1)\setminus X)$};
\draw (\D+\d*3/4,-\h) node (A1') { };
\draw (\L/2 + \D/2 + \d/2,-\h) node (B1') { } ;
\draw[->,>=latex] (A) to[out=-90, in =90,distance=15] (A1');
\draw[dotted,->,>=latex] (B1) to[out=-90, in =90, distance=15] (B1');
\draw[dotted,->,>=latex] (B2) to[out=-90, in =90, distance = 10] (B2');
\draw[->,>=latex] (C) to[out=-90, in =90] (A2');
\end{footnotesize}

\draw [->] (-0.05/2,0) to [out = -135,in =135] (-0.05/2,-\h);
\draw (-0.2,-\h/2) node {$g$}; 
\draw (\L+0.2,-\h/2) node {}; 

\end{tikzpicture}
\caption{Example of AIETs $f$ and $g$ spanning a free group of rank 2} \label{fig_AIET}
\end{figure}

\end{proof}

\bibliographystyle{abbrv}

\bibliography{references}

\begin{thebibliography}{10}

\bibitem{KatokQuestion}
In the conference report, the question is not clearly attributed to someone,
  but K. Fr\k{a}czek, who was there, confirmed to the author that it was A.
  Katok who asked the question. Of course, other mathematicians might have
  posed this question before but it was not investigated.

\bibitem{ArnouxFIET}
P.~Arnoux.
\newblock Un invariant pour les échanges d'intervalles et les flots sur les
  surfaces. {T}hèse. {U}niversité de {R}eims, 1981.

\bibitem{Arnoux_eif}
P.~Arnoux.
\newblock Échanges d'intervalles et flots sur les surfaces.
\newblock {\em Theorie ergodique, {Semin}. {Les} {Plans}-sur-{Bex} 1980,
  {Monogr}. {L}'{Enseign}. {Math}. (French)}, 29:5--38, 1981.

\bibitem{Arnoux}
P.~Arnoux.
\newblock Un exemple de semi-conjugaison entre un \'echange d'intervalles et
  une translation sur le tore.
\newblock {\em Bulletin de la Soci\'et\'e Math\'ematique de France},
  116(4):489--500, 1988.

\bibitem{AY}
P.~{Arnoux} and J.-C. {Yoccoz}.
\newblock {Construction de diff\'eomorphismes pseudo-Anosov}.
\newblock {\em {C. R. Acad. Sci., Paris, S\'er. I}}, 292:75--78, 1981.

\bibitem{Boshernitzan}
M.~Boshernitzan.
\newblock Subgroup of interval exchanges generated by torsion elements and
  rotations.
\newblock {\em Proc. Am. Math. Soc.}, 144(6):2565--2573, 2016.

\bibitem{Coffrey}
J.~Coffey.
\newblock Some remarks concerning an example of a minimal, non-uniquely ergodic
  interval exchange transformation.
\newblock {\em Math. Z.}, 199(4):577--580, 1988.

\bibitem{DFG}
F.~Dahmani, K.~Fujiwara, and V.~Guirardel.
\newblock Free groups of interval exchange transformations are rare.
\newblock {\em {Groups Geom. Dyn.}}, 7(4):883--910, 2013.

\bibitem{DeLaHarpe}
P.~de~la Harpe.
\newblock {\em Topics in geometric group theory}.
\newblock Chicago Lectures in Mathematics. University of Chicago Press,
  Chicago, IL, 2000.

\bibitem{KatokStepin}
A.~B. Katok and A.~M. Stepin.
\newblock Approximations in ergodic theory.
\newblock {\em Usp. Mat. Nauk}, 22(5(137)):81--106, 1967.

\bibitem{Keane75}
M.~Keane.
\newblock Interval exchange transformations.
\newblock {\em Mathematische Zeitschrift}, 141:25--32, 1975.

\bibitem{Keane77}
M.~Keane.
\newblock Non-ergodic interval exchange transformations.
\newblock {\em Isr. J. Math.}, 26:188--196, 1977.

\bibitem{KeynesNewton}
H.~B. Keynes and D.~Newton.
\newblock A 'minimal', non-uniquely ergodic interval exchange transformation.
\newblock {\em Math. Z.}, 148:101--105, 1976.

\bibitem{Sah}
C.-H. Sah.
\newblock Scissors congruences of the interval. {P}reprint, 1981.

\bibitem{Sataev}
E.~A. Sataev.
\newblock On the number of invariant measures for flows on orientable surfaces.
\newblock {\em Math. USSR, Izv.}, 9:813--830, 1976.

\bibitem{Veech}
W.~A. Veech.
\newblock Strict ergodicity in zero dimensional dynamical systems and the
  {Kronecker}-{Weyl} theorem {{\(\bmod 2\)}}.
\newblock {\em Trans. Am. Math. Soc.}, 140:1--33, 1969.

\bibitem{Viana}
M.~Viana.
\newblock Ergodic theory of interval exchange maps.
\newblock {\em Revista Matem\'atica de la Universidad Complutense de Madrid},
  19:7--100, 2006.

\bibitem{Vorobets}
Y.~Vorobets.
\newblock On the commutator group of the group of interval exchange
  transformations.
\newblock {\em Proc. Steklov Inst. Math.}, 297:285--296, 2017.

\end{thebibliography}

\end{document}